\documentclass{article}%
\usepackage{amssymb}
\usepackage{amsmath}
\usepackage{amsfonts}
\usepackage{graphicx}%
\providecommand{\U}[1]{\protect\rule{.1in}{.1in}}
\newtheorem{theorem}{Theorem}

\newtheorem{definition}[theorem]{Definition}

\newtheorem{remark}[theorem]{Remark}

\begin{document}

\title{An abstract theorem on the existence of hylomorphic solitons}
\author{Vieri Benci$^{\ast}$, Donato Fortunato$^{\ast\ast}$\\$^{\ast}$Dipartimento di Matematica Applicata \textquotedblleft U.
Dini\textquotedblright\\Universit\`{a} degli Studi di Pisa \\Via Filippo Buonarroti 1/c, 56127 Pisa, Italy\\e-mail: benci@dma.unipi.it\\College of Science, Department of Mathematics \\King Saud University \\Riyadh, 11451, Saudi Arabia \\$^{\ast\ast}$Dipartimento di Matematica \\Universit\`{a} degli Studi di Bari Aldo Moro\\Via Orabona 4, 70125 Bari, Italy\\e-mail: fortunat@dm.uniba.it}
\maketitle

\begin{abstract}
In this paper we prove an abstract theorem which can be used to study the
existence of solitons for various dynamical systems described by partial
differential equations. We also give an idea of how the abstract theorem can
be applied to prove the existence of solitons in some dynamical systems.

\end{abstract}
\tableofcontents


\bigskip


\section{Introduction}

In some recent papers ( \cite{BBBM},\ \cite{BBGM}, \cite{milano},
\cite{befobeam}, \cite{befolat}, \cite{befoQ}) the existence of solitons has
been proved using very similar techniques. In this paper we prove a general,
abstract theorem which applies to most of the situations analyzed in the
mentioned papers.

In section \ref{be} we give an abstract definition of soliton. Section 3 is
devoted to the main result, namely an existence result for hylomorphic
solitons (Theorem \ref{astra1}). These solitons are minimizers, satisfying
suitable stability properties, of a constrained functional. The proof of
Theorem \ref{astra1} is carried out in two steps: in the first step the
research of the minimizers of a constrained functional is reduced to the study
of the minimizers of a suitable free functional (Theorem \ref{furbo}). In the
second step stability properties of these minimizers are proved (subsection 3.2).

In section \ref{application} we give an idea (the complete proofs being
contained in the quoted papers) of how the abstract theorem \ref{astra1} can
be used to prove the existence of solitons in some dynamical systems, namely
for the nonlinear Schr\"{o}dinger equation, the nonlinear wave equation, the
nonlinear beam equation and the Klein-Gordon-Maxwell equation. In a
forthcoming paper we shall use theorem \ref{astra1} to study new situations.

\section{An abstract definition of solitons\label{be}}

Solitons are particular \textit{states} of a dynamical system described by one
or more partial differential equations. Thus, we assume that the states of
this system are described by one or more \textit{fields} which mathematically
are represented by functions
\begin{equation}
\mathbf{u}:\mathbb{R}^{N}\rightarrow V \label{lilla}%
\end{equation}
where $V$ is a vector space with norm $\left\vert \ \cdot\ \right\vert _{V}$
which is called the internal parameters space. We assume the system to be
deterministic; this means that it can be described as a dynamical system
$\left(  X,\gamma\right)  $ where $X$ is the set of the states and
$\gamma:\mathbb{R}\times X\rightarrow X$ is the time evolution map. If
$\mathbf{u}_{0}(x)\in X,$ the evolution of the system will be described by the
function
\begin{equation}
\mathbf{u}\left(  t,x\right)  :=\gamma_{t}\mathbf{u}_{0}(x). \label{flusso}%
\end{equation}
We assume that the states of $X$ have "finite energy" so that they decay at
$\infty$ sufficiently fast. Roughly speaking, the solitons are "bump"
solutions characterized by some form of stability.

To define them at this level of abstractness, we need to recall some well
known notions in the theory of dynamical systems.

\begin{definition}
A set $\Gamma\subset X$ is called \textit{invariant} if $\forall\mathbf{u}%
\in\Gamma,\forall t\in\mathbb{R},\ \gamma_{t}\mathbf{u}\in\Gamma.$
\end{definition}

\begin{definition}
Let $\left(  X,d\right)  $ be a metric space and let $\left(  X,\gamma\right)
$ be a dynamical system. An invariant set $\Gamma\subset X$ is called stable,
if $\forall\varepsilon>0,$ $\exists\delta>0,\;\forall\mathbf{u}\in X$,
\[
d(\mathbf{u},\Gamma)\leq\delta,
\]
implies that
\[
\forall t\geq0,\text{ }d(\gamma_{t}\mathbf{u,}\Gamma)\leq\varepsilon.
\]

\end{definition}

\bigskip

For every $z\in\mathbb{Z}^{N},$ and $\mathbf{u}\in X$, we set
\begin{equation}
\left(  g_{z}\mathbf{u}\right)  \left(  x\right)  =\mathbf{u}\left(
x-Az\right)  . \label{ggg}%
\end{equation}
where $A$ is an inverible matrix; such a $g_{z}$ will be called
\textit{lattice transformation}. We set
\begin{equation}
G=\left\{  g_{z}|\ z\in\mathcal{G}\right\}  ; \label{gg}%
\end{equation}
where $\mathcal{G}$ is a subgroup of $\left(  \mathbb{Z}^{N},+\right)  .\ G$
is a group of transformations acting on the space $X;$ actually it is a linear
representation of the group $\mathcal{G}.$

\begin{definition}
A non-empty subset $\Gamma\subset X$ is called $G$-invariant if
\[
\forall\mathbf{u}\in\Gamma,\ \forall z\in\mathbb{Z}^{N},\ g_{z}\mathbf{u}%
\in\Gamma.
\]

\end{definition}

\begin{definition}
\label{gcompatto}A closed $G$-invariant set $\Gamma\subset X$ is called
$G$-compact if for any sequence $\mathbf{u}_{n}(x)$ in $\Gamma$ there is a
sequence $g_{n}\in G,$ such that $\mathbf{u}_{n}(g_{n}x)$ has a converging subsequence.
\end{definition}

Now we are ready to give the definition of soliton:

\begin{definition}
\label{ds} A state $\mathbf{u}(x)\in X$ is called soliton if there is an
invariant set $\Gamma$ such that

\begin{itemize}
\item (i) $\forall t,\ \gamma_{t}\mathbf{u}(x)\in\Gamma,$

\item (ii) $\Gamma$ is stable,

\item (iii) $\Gamma$ is $G$-compact.
\end{itemize}
\end{definition}

\begin{remark}
The above definition needs some explanation. For simplicity, we assume that
$\Gamma$ is a manifold (actually, in many concrete models, this is the generic
case). Then (iii) implies that $\Gamma$ is finite dimensional. Since $\Gamma$
is invariant, $\mathbf{u}_{0}\in\Gamma\Rightarrow\gamma_{t}\mathbf{u}_{0}%
\in\Gamma$ for every time. Thus, since $\Gamma$ is finite dimensional, the
evolution of $\mathbf{u}_{0}$ is described by a finite number of parameters$.$
The dynamical system $\left(  \Gamma,\gamma\right)  $\ behaves as a point in a
finite dimensional phase space. By the stability of $\Gamma$, a small
perturbation of $\mathbf{u}_{0}$ remains close to $\Gamma.$ However, in this
case, its evolution depends on an infinite number of parameters. Thus, this
system appears as a finite dimensional system with a small perturbation.
\end{remark}

\bigskip

\section{Existence of hylomorphic solitons}

\bigskip

We now assume that the dynamical system $\left(  X,\gamma\right)  $ has $l+1$
integrals of motion. One of them will be called energy and it will be denoted
by $E$; the set of the other integrals can be considered as a functional
\[
C:X\rightarrow\mathbb{R}^{l}%
\]
and it will be called hylenic charge. At this level of abstractness, the name
energy and hylenic charge are conventional, but $E$ and $C$ need to satisfy
different assumptions see assumption (EC-3) below. In our applications to
PDE,\ $E$ will be the usual energy. The name hylenic charge has been
introduced in \cite{milano}, \cite{BBBM}, \cite{hylo} and here this name will
denote just a set of other integrals.

The presence of $E$ and $C$ allows to give the following definition of
hylomorphic soliton.

\begin{definition}
\label{tdc}A soliton $\mathbf{u}_{0}\in X$ is hylomorphic if $\Gamma$ (as in
Def. \ref{ds}) has the following structure%
\[
\Gamma=\Gamma\left(  e_{0},c_{0}\right)  =\left\{  \mathbf{u}\in
X\ |\ E(\mathbf{u})=e_{0},\ C(\mathbf{u})=c_{0}\right\}
\]
where%
\begin{equation}
e_{0}=\min\left\{  E(\mathbf{u})\ |\ C(\mathbf{u})=c_{0}\right\}  \label{min}%
\end{equation}
for some $c_{0}\in\mathbb{R}^{l}$.
\end{definition}

Notice that, by (\ref{min}), we have that a hylomorphic soliton $\mathbf{u}%
_{0}$ satisfies the following nonlinear eigenvalue problem:%
\[
E^{\prime}(\mathbf{u}_{0})=\sum_{j=1}^{l}\lambda_{j}C_{j}^{\prime}%
(\mathbf{u}_{0})
\]
where we have set $C(\mathbf{u})=(C_{1}(\mathbf{u}),...,C_{l}(\mathbf{u}))$

Clearly, for a given $c_{0},$ the minimum of $E$ might not exist; moreover,
even if the minimum exists, it is possible that $\Gamma\ $ does not satisfy
(ii) or (iii) of def. \ref{ds}.

In order to prove an existence result for hylomorphic solitons, we impose some
assumptions to $E$ and $C$; to do this we need some other definitions:

\begin{definition}
We say that a functional $F$ on $X$ has the splitting property if given a
sequence $\mathbf{u}_{n}=\mathbf{u}+\mathbf{w}_{n}\in X$ such that
$\mathbf{w}_{n}$ converges weakly to $0$, we have that
\begin{equation}
F(\mathbf{u}_{n})=F(\mathbf{u})+F(\mathbf{w}_{n})+o(1)
\end{equation}

\end{definition}

\bigskip

\begin{remark}
Every quadratic form satisfies the splitting property; in fact, in this case,
we have that $F(\mathbf{u}):=\left\langle L\mathbf{u},\mathbf{u}\right\rangle
$ for some continuous selfajoint operator $L;$ then, given a sequence
$\mathbf{u}_{n}=\mathbf{u}+\mathbf{w}_{n}$ with $\mathbf{w}_{n}\rightharpoonup
0$ weakly, we have that%
\begin{align*}
F(\mathbf{u}_{n})  &  =\left\langle L\mathbf{u},\mathbf{u}\right\rangle
+\left\langle L\mathbf{w}_{n},\mathbf{w}_{n}\right\rangle +2\left\langle
L\mathbf{u},\mathbf{w}_{n}\right\rangle \\
&  =F(\mathbf{u})+F(\mathbf{w}_{n})+o(1)
\end{align*}

\end{remark}

\bigskip

Now we can formulate the properties we require:

\begin{itemize}
\item \textit{(EC-1) there are }$l+1$\textit{\ prime integrals }%
$E(\mathbf{u})$\textit{\ and }$C(\mathbf{u})=(C_{1}(\mathbf{u}),...,C_{l}%
(\mathbf{u}))$\textit{\ of the dynamical system }$(X,\gamma)$\textit{\ such
that}%
\[
E(0)=0,\ C(0)=0;\ E^{\prime}(0)=0;\ C^{\prime}(0)=0
\]

\item \textit{(EC-2) }$E(u)$\textit{\ and }$C(u)$\textit{\ are }%
$G$\textit{-invariant.}

\item \textit{(EC-3)(coercivity assumption) suppose that there exists }%
$a\geq0,$ and $s\geq1$ such that

\begin{itemize}
\item (i) $E(\mathbf{u})+a\left\vert C(\mathbf{u})\right\vert ^{s}\geq0;$

\item (ii) \textit{if }$\left\Vert \mathbf{u}\right\Vert \rightarrow\infty
,\ $\textit{then} $E(\mathbf{u})+a\left\vert C(\mathbf{u})\right\vert
^{s}\rightarrow\infty$

\item (iii) \textit{for any} \textit{sequence }$\mathbf{u}_{n}$ \textit{in
}$X$ \textit{such that} $E(\mathbf{u}_{n})+a\left\vert C(\mathbf{u}%
_{n})\right\vert ^{s}\rightarrow0,\ $\textit{we have that }$\mathbf{u}%
_{n}\rightarrow0$
\end{itemize}

\item \textit{(EC-4) }$E$\textit{\ and }$|C|$\textit{\ satisfy the splitting
property.}
\end{itemize}

\bigskip

\begin{definition}
\label{na}A norm (or a seminorm) $\left\Vert \cdot\right\Vert _{\sharp}$ is
called a auxiliary norm (or seminorm) on $X$ if it satisfies the following property:

\begin{itemize}
\item given any sequence $\left\Vert \mathbf{u}_{n}\right\Vert $ bounded in
$X$ such that $\left\Vert \mathbf{u}_{n}\right\Vert _{\sharp}\geq\delta>0,$ we
can extract a subsequence $\mathbf{u}_{n_{k}}$ and we can take a sequence
$g_{k}\in G$ such that $g_{k}\mathbf{u}_{n_{k}}$ is weakly convergent to some
$\mathbf{\bar{u}}\neq0$.
\end{itemize}
\end{definition}

\bigskip

The notion of auxiliary norm is related to a result obtained by Lieb in
\cite{lieb}. In our abstract scheme, this norm allows to define the following
number:%
\begin{equation}
\Lambda_{0}:=\ \underset{\left\Vert \mathbf{u}\right\Vert _{\sharp}%
\rightarrow0}{\lim\inf}\frac{E(\mathbf{u})}{\left\vert C(\mathbf{u}%
)\right\vert }. \label{hylo}%
\end{equation}
Now we can state the main result:\bigskip

\begin{theorem}
\label{astra1} Assume that the dynamical system $(X,\gamma)$ satisfies
(EC-1),...,(EC-4). Moreover assume that%
\begin{equation}
\underset{\mathbf{u}\in X}{\inf}\frac{E(\mathbf{u})}{\left\vert C(\mathbf{u}%
)\right\vert }<\Lambda_{0}. \label{hh}%
\end{equation}
\newline Then, the dynamical system $(X,\gamma)$ admits a continuous family
$\mathbf{u}_{\delta}$ ($\delta\in\left(  0,\bar{\delta}\right)  ,$
$\bar{\delta}>0$) of independent, hylomorphic solitons (two solitons
$\mathbf{u}_{\delta_{1}},\mathbf{u}_{\delta_{2}}$ are called independent if
$\mathbf{u}_{\delta_{1}}\neq g\mathbf{u}_{\delta_{2}}$ for every $g\in G$).
\end{theorem}

\bigskip

We will prove Theorem \ref{astra1} in two steps:

\begin{itemize}
\item first step: existence of a set of minimizers $\Gamma$ as in def.
\ref{tdc}

\item second step: proof of the stability of $\Gamma$
\end{itemize}

\bigskip

\subsection{An existence result for constrained minimizers}

In this section, we will be concerned with the first step. We need some
definitions. These definitions are related to a couple $(X,G)$ where $G$ is a
group acting on $X.$ In our applications $G$ will be the group (\ref{ggg}),
but these results hold also in an abstract framework.

\begin{definition}
Let $G$ be a group acting on $X.$ A sequence $\mathbf{u}_{n}$ in $X$ is called
$G$\emph{-compact }if we can extract a subsequence $\mathbf{u}_{n_{k}} $ such
that there exists a sequence $g_{k}\in G$ such that $g_{k}\mathbf{u}_{n_{k}}$
is convergent.
\end{definition}

\begin{definition}
A functional $J$ on $X$ is called $G$\emph{-invariant} if
\[
\forall g\in G,\text{ }\forall\mathbf{u}\in X,\ J\left(  g\mathbf{u}\right)
=J\left(  \mathbf{u}\right)  .
\]

\end{definition}

\begin{definition}
A functional $J$ on $X$ is called $G$\emph{-compact} if any minimizing
sequence $\mathbf{u}_{n}$ is\emph{\ }$G$-compact.
\end{definition}

\bigskip

Notice that a $G$-compact functional has a $G$-compact set of minimizers (see
def. \ref{gcompatto}).

Now, we will prove the following existence result:

\begin{theorem}
\label{furbo} Assume that $E$ and $C$ satisfy (EC-1),...,(EC-4) and
(\ref{hh}). Then there are $\bar{\delta}>0$ and a family of values of the
charge $c_{\delta},\ \delta\in\left(  0,\bar{\delta}\right)  ,$ such that, for
all $\delta\in\left(  0,\bar{\delta}\right)  ,$ the minimum
\begin{equation}
e_{\delta}=\min\left\{  E(\mathbf{u})\ |\ C(\mathbf{u})=c_{\delta}\right\}
\end{equation}
exists and the set of minimizers $\Gamma_{\delta}$ is $G$-compact. Moreover,
$\Gamma_{\delta}\ $can be characterized also as the set of minimizers of the
$G$-compact functional%
\[
J_{\delta}(\mathbf{u})=\Lambda\left(  \mathbf{u}\right)  +\delta
\Phi(\mathbf{u})
\]
where%
\[
\Lambda\left(  \mathbf{u}\right)  =\frac{E(\mathbf{u})}{\left\vert
C(\mathbf{u})\right\vert }%
\]%
\[
\Phi(\mathbf{u})=E(\mathbf{u})+2a\left\vert C(\mathbf{u})\right\vert ^{s}.
\]

\end{theorem}

\textbf{Proof.} Now take the functional
\[
J_{\delta}(\mathbf{u})=\Lambda\left(  \mathbf{u}\right)  +\delta
\Phi(\mathbf{u})
\]
where $\delta\in\left(  0,\bar{\delta}\right)  $ is so small that
\begin{equation}
\underset{\mathbf{u}\in X}{\inf}J_{\delta}\left(  \mathbf{u}\right)
<\Lambda_{0}. \label{fava}%
\end{equation}
Such $\bar{\delta}>0$ does exist by virtue of assumption (\ref{hh}).

It is easy to show that
\begin{equation}
J_{\delta}\left(  \mathbf{u}\right)  \geq\frac{\delta}{2}\Phi(\mathbf{u})-M
\label{pli}%
\end{equation}
where $M$ is a suitable constant; in fact by (EC-3)(i) we have that%
\begin{equation}
E(\mathbf{u})\geq-a\left\vert C(\mathbf{u})\right\vert ^{s} \label{aaa}%
\end{equation}
and hence%
\begin{equation}
\frac{E(\mathbf{u})}{\left\vert C(\mathbf{u})\right\vert }\geq-a\left\vert
C(\mathbf{u})\right\vert ^{s-1} \label{ccc}%
\end{equation}

Then, by (\ref{aaa}) and (\ref{bbb})
\begin{align*}
J_{\delta}(\mathbf{u})  &  =\frac{E(\mathbf{u})}{\left\vert C(\mathbf{u}%
)\right\vert }+\delta\Phi(\mathbf{u})\geq-a\left\vert C(\mathbf{u})\right\vert
^{s-1}+\frac{\delta}{2}\left[  E(\mathbf{u})+2a\left\vert C(\mathbf{u}%
)\right\vert ^{s}\right]  +\frac{\delta}{2}\Phi(\mathbf{u})\\
&  \geq-a\left\vert C(\mathbf{u})\right\vert ^{s-1}+\frac{\delta}{2}\left[
-a\left\vert C(\mathbf{u})\right\vert ^{s}+2a\left\vert C(\mathbf{u}%
)\right\vert ^{s}\right]  +\frac{\delta}{2}\Phi(\mathbf{u})\\
&  =-a\left\vert C(\mathbf{u})\right\vert ^{s-1}+\frac{a\delta}{2}\left\vert
C(\mathbf{u})\right\vert ^{s}+\frac{\delta}{2}\Phi(\mathbf{u})\geq\frac
{\delta}{2}\Phi(\mathbf{u})-M
\end{align*}
where%
\[
M=\ -a\underset{t\geq0}{\min}\left(  \frac{\delta}{2}t^{s}-t^{s-1}\right)  .
\]

Next, we will prove that $J_{\delta}$ is $G$-compact.\ Let $\mathbf{u}_{n}$ be
a minimizing sequence of $J_{\delta}.$ We have to prove that there exists a
sequence $g_{k}\in G$ and a subsequence $\mathbf{u}_{n_{k}}$ such that
$\mathbf{u}_{k}^{\prime}:=g_{k}\mathbf{u}_{n_{k}}$ is convergent.

By (\ref{fava}), there exists $\eta>0$ such that, for $n$ sufficiently large,
\[
\frac{E(\mathbf{u}_{n})}{\left\vert C(\mathbf{u}_{n})\right\vert }+\delta
\Phi(\mathbf{u}_{n})<\Lambda_{0}-\eta.
\]
So we have that for $n$ sufficiently large
\[
\frac{E(\mathbf{u}_{n})}{\left\vert C(\mathbf{u}_{n})\right\vert }<\Lambda
_{0}-\eta
\]
and hence, by the definition of $\Lambda_{0}$, we have that
\begin{equation}
\left\Vert \mathbf{u}_{n}\right\Vert _{\sharp}\geq b \label{bbb}%
\end{equation}
for some $b>0.\ $By (EC-3)(ii) and (\ref{pli}), we have that $\left\Vert
\mathbf{u}_{n}\right\Vert $ is bounded. Then, by Def \ref{na}, we can extract
a subsequence $\mathbf{u}_{n_{k}}$ and we can take a sequence $g_{k}\in G$
such that $\mathbf{u}_{k}^{\prime}:=g_{k}\mathbf{u}_{n_{k}}$ is weakly
convergent to some
\begin{equation}
\mathbf{\bar{u}}\neq0. \label{agg}%
\end{equation}
We can write
\[
\mathbf{u}_{n}^{\prime}=\mathbf{\bar{u}}+\mathbf{w}_{n}%
\]
with $\mathbf{w}_{n}\rightharpoonup0$ weakly. We want to prove that
$\mathbf{w}_{n}\rightarrow0$ strongly. First of all we will show that%
\[
\lim\Phi(\mathbf{\bar{u}}+\mathbf{w}_{n})\geq\Phi(\mathbf{\bar{u}})+\lim
\Phi(\mathbf{w}_{n}).
\]
In fact, using (EC-4) we have that
\begin{align}
\lim\Phi(\mathbf{\bar{u}}+\mathbf{w}_{n})  &  =\lim\left(  E(\mathbf{\bar{u}%
}+\mathbf{w}_{n})+2a\left\vert C(\mathbf{\bar{u}}+\mathbf{w}_{n})\right\vert
^{s}\right) \nonumber\\
&  =E(\mathbf{\bar{u}})+\lim E(\mathbf{w}_{n})+2a\lim\left(  \left\vert
C(\mathbf{\bar{u}})\right\vert +\left\vert C(\mathbf{w}_{n})\right\vert
\right)  ^{s}\nonumber\\
&  \geq E(\mathbf{\bar{u}})+\lim E(\mathbf{w}_{n})+2a\lim\left(  \left\vert
C(\mathbf{\bar{u}})\right\vert ^{s}+\left\vert C(\mathbf{w}_{n})\right\vert
^{s}\right) \nonumber\\
&  =E(\mathbf{\bar{u}})+2a\left\vert C(\mathbf{\bar{u}})\right\vert ^{s}+\lim
E(\mathbf{w}_{n})+2a\lim\left\vert C(\mathbf{w}_{n})\right\vert ^{s}%
\nonumber\\
&  =\Phi(\mathbf{\bar{u}})+\lim\Phi(\mathbf{w}_{n}). \label{qw}%
\end{align}

Next we show that%
\begin{equation}
C(\mathbf{\bar{u}+\mathbf{w}_{n})}\text{ does not converge to }0. \label{not}%
\end{equation}
Arguing by contradiction assume that $C(\mathbf{\bar{u}+\mathbf{w}_{n})}$
converges to $0.$ Then, since $\mathbf{\bar{u}}+\mathbf{\mathbf{w}_{n}}$ is a
minimizing sequence for $J_{\delta},$ also $E(\mathbf{\bar{u}+\mathbf{w}_{n}%
)}$\ and $\Phi(\mathbf{\bar{u}+\mathbf{w}_{n})}$ converge to $0$. Then, by
(EC-3)(iii), we get
\begin{equation}
\mathbf{\bar{u}+\mathbf{w}_{n}}\rightarrow0\text{ in }X \label{s}%
\end{equation}
From (\ref{s}) and since $\mathbf{\mathbf{w}_{n}}\rightarrow0$ weakly in $X, $
we have that $\mathbf{\bar{u}}=0$, contradicting (\ref{agg}). Then, by
(\ref{not}),
\begin{equation}
C(\mathbf{\bar{u}+\mathbf{w}_{n})}=C(\mathbf{\bar{u}})+C(\mathbf{w}%
_{n})+o(1)\geq const.>0 \label{fir}%
\end{equation}
up to a subsequence. Now we set%
\begin{align*}
j_{\delta}  &  =\lim J_{\delta}(\mathbf{u}_{n}^{\prime});\ e_{\delta
}=E(\mathbf{\bar{u}});\ c_{\delta}=\left\vert C(\mathbf{\bar{u}})\right\vert
\\
e_{1}  &  =\lim E(\mathbf{w}_{n});\ c_{1}=\lim\left\vert C(\mathbf{w}%
_{n})\right\vert .
\end{align*}

By (\ref{qw}) and (\ref{fir}), we have that%
\begin{align}
j_{\delta}  &  =\lim\left[  \frac{E(\mathbf{u}_{n}^{\prime})}{\left\vert
C(\mathbf{u}_{n}^{\prime})\right\vert }+\delta\Phi(\mathbf{u}_{n}^{\prime
})\right] \nonumber\\
&  =\lim\frac{E(\mathbf{\bar{u}})+E(\mathbf{w}_{n})+o(1)}{\left\vert
C(\mathbf{\bar{u}})+C(\mathbf{w}_{n})+o(1)\right\vert }+\delta\lim
\Phi(\mathbf{\bar{u}}+\mathbf{w}_{n})\nonumber\\
&  \geq\frac{\ e_{\delta}+e_{1}}{c_{\delta}+c_{1}}+\delta\lim\Phi
(\mathbf{w}_{n})\ +\delta\Phi(\mathbf{\bar{u}}). \label{pipa}%
\end{align}
Now we want to prove that
\begin{equation}
\frac{\ e_{1}}{c_{1}}\geq\frac{\ e_{\delta}}{c_{\delta}}. \label{impo2}%
\end{equation}

We argue indirectly and we suppose that
\begin{equation}
\frac{\ e_{\delta}}{c_{\delta}}>\frac{\ e_{1}}{c_{1}}. \label{impo}%
\end{equation}
By the above inequality it follows that
\begin{equation}
\frac{\ e_{\delta}+e_{1}}{c_{\delta}+c_{1}}=\frac{\frac{e_{\delta}}{c_{\delta
}}\ c_{\delta}+\frac{e_{1}}{c_{1}}\ c_{1}}{c_{\delta}+c_{1}}>\frac{\frac
{e_{1}}{c_{1}}\ c_{\delta}+\frac{e_{1}}{c_{1}}\ c_{1}}{c_{\delta}+c_{1}}%
=\frac{e_{1}}{c_{1}} \label{ava}%
\end{equation}
and hence%
\begin{align*}
j_{\delta}  &  \geq\frac{\ e_{\delta}+e_{1}}{c_{\delta}+c_{1}}+\delta\lim
\Phi(\mathbf{w}_{n})\ +\delta\Phi(\mathbf{\bar{u}})\\
&  >\frac{e_{1}}{c_{1}}+\delta\lim\Phi(\mathbf{w}_{n})+\delta\Phi
(\mathbf{\bar{u}})\\
&  =\lim J_{\delta}(\mathbf{w}_{n})+\delta\Phi(\mathbf{\bar{u}})\geq
j_{\delta}+\delta\Phi(\mathbf{\bar{u}})>j_{\delta}.
\end{align*}
So (\ref{impo}) cannot occur and then we have (\ref{impo2}). In this case,
arguing as in (\ref{ava})%
\begin{equation}
\frac{\ e_{\delta}+e_{1}}{c_{\delta}+c_{1}}\geq\frac{e_{\delta}}{c_{\delta}}%
\end{equation}
ans so, using (\ref{pipa}) and the above inequality,
\begin{align*}
j_{\delta}  &  \geq\frac{\ e_{\delta}+e_{1}}{c_{\delta}+c_{1}}+\delta
\Phi(\mathbf{\bar{u}})+\delta\lim\Phi(\mathbf{w}_{n})\\
&  \geq\frac{e_{\delta}}{c_{\delta}}+\delta\Phi(\mathbf{\bar{u}})+\delta
\lim\Phi(\mathbf{w}_{n})\\
&  =J_{\delta}(\mathbf{\bar{u}})+\delta\lim\Phi(\mathbf{w}_{n})\geq j_{\delta
}+\delta\lim\Phi(\mathbf{w}_{n}).
\end{align*}
Then%
\[
\delta\lim\Phi(\mathbf{w}_{n})\leq0
\]
and by (EC-3)(iii) $\mathbf{w}_{n}\rightarrow0\ $and hence $\mathbf{u}%
_{n}^{\prime}\rightarrow\mathbf{\bar{u}}$ strongly. Thus $J_{\delta}$ is
$G$-compact and the set of minimizer $\Gamma_{\delta}$ is not empty. Clearly,
if $\mathbf{u}\in\Gamma_{\delta},$ it turns out that $\mathbf{u}$ minimises
also the functional%
\[
\frac{E(\mathbf{u})}{c_{\delta}}+\delta\left[  E(\mathbf{u})+ac_{\delta}%
^{s}\right]  =\left(  \frac{1}{c_{\delta}}+\delta\right)  E(\mathbf{u})+\delta
ac_{\delta}^{s}%
\]
on the set $\left\{  \mathbf{u}\in X\ |\ C(\mathbf{u})=c_{\delta}\right\}  $
and hence it minimizes $E(\mathbf{u})$ on this set.

$\square$

In the next section we consider the second step, namely we prove that
$\Gamma_{\delta}$ is stable and then we complete the proof of Theorem
\ref{astra1}

\subsection{A stability result}

We need the (well known) Liapunov theorem in following form:

\begin{theorem}
\label{propV}Let $\Gamma$ be an invariant set and assume that there exists a
differentiable function $V$ (called a Liapunov function) such that

\begin{itemize}
\item (a) $V(\mathbf{u})\geq0$ and\ $V(\mathbf{u})=0\Leftrightarrow u\in
\Gamma$

\item (b) $\partial_{t}V(\gamma_{t}\left(  \mathbf{u}\right)  )\leq0$

\item (c) $V(\mathbf{u}_{n})\rightarrow0\Leftrightarrow d(\mathbf{u}%
_{n},\Gamma)\rightarrow0.$
\end{itemize}

\noindent Then $\Gamma$ is stable.
\end{theorem}

\textbf{Proof. }For completeness, we give a proof of this well known result.
Arguing by contradiction, assume that $\Gamma,$ satisfying the assumptions of
Th. \ref{propV}, is not stable. Then there exists $\varepsilon>0$ and
sequences $\mathbf{u}_{n}\in X$ and $t_{n}>0$ such that
\begin{equation}
d(\mathbf{u}_{n},\Gamma)\rightarrow0\text{ and }d(\gamma_{t_{n}}\left(
\mathbf{u}_{n}\right)  ,\Gamma)>\varepsilon. \label{bingo}%
\end{equation}
Then we have%
\[
d(\mathbf{u}_{n},\Gamma)\rightarrow0\Longrightarrow V(\mathbf{u}%
_{n})\rightarrow0\Longrightarrow V(\gamma_{t_{n}}\left(  \mathbf{u}%
_{n}\right)  )\rightarrow0\Longrightarrow d(\gamma_{t_{n}}\left(
\mathbf{u}_{n}\right)  ,\Gamma)\rightarrow0
\]
where the first and the third implications are consequence of property (c).
The second implication follows from property (b). Clearly, this fact
contradicts (\ref{bingo}).

$\square$

The following Theorem holds:

\begin{theorem}
\label{astraco}Assume (EC-1) and (EC-2). For $\mathbf{u}\in X$ and
$e_{0},c_{0}\in$ $\mathbb{R}$, we set%
\begin{equation}
V\left(  \mathbf{u}\right)  =\left(  E\left(  \mathbf{u}\right)
-e_{0}\right)  ^{2}+\left(  C\left(  \mathbf{u}\right)  -c_{0}\right)  ^{2}.
\end{equation}
If $V$ is$\ G$-compact and
\begin{equation}
\Gamma=\left\{  \mathbf{u}\in X:E(\mathbf{u})=e_{0},\ C(\mathbf{u}%
)=c_{0}\right\}  \neq\varnothing, \label{gamma}%
\end{equation}
then every $\mathbf{u}\in\Gamma$ is a soliton.
\end{theorem}

\textbf{Proof}: We have to prove that $\Gamma$ in (\ref{gamma}) satisfies
(i),(ii) and (iii) of Def. (\ref{ds}). The property (iii), namely the fact
that $\Gamma$ is G-compact, is a trivial consequence of the fact that $\Gamma$
is the set of minimizers of a G-compact functional $V.$ The invariance
property (i) is clearly satisfied since $E$ and $C$ are constants of the
motion. It remains to prove (ii), namely that $\Gamma$ is stable. To this end
we shall use Th. \ref{propV}. So we need to show that $V(\mathbf{u}) $
satisfies (a), (b) and (c). Statements (a) and (b) are trivial. Now we prove
(c). First we show the implication $\Rightarrow.$ Let $\mathbf{u}_{n}$ be a
sequence such that $V(\mathbf{u}_{n})\rightarrow0.$ By contradiction we assume
that $d(\mathbf{u}_{n},\Gamma)\nrightarrow0,$ namely that there is a
subsequence $\mathbf{u}_{n}^{^{\prime}}$ such that
\begin{equation}
d(\mathbf{u}_{n}^{\prime},\Gamma)\geq a>0. \label{kaka}%
\end{equation}
Since $V(\mathbf{u}_{n})\rightarrow0$ also $V(\mathbf{u}_{n}^{\prime
})\rightarrow0,$ and, since $V$ is $G$ compact, there exists a sequence
$g_{n}$ in $G$ such that, for a subsequence $\mathbf{u}_{n}^{\prime\prime}$,
we have $g_{n}\mathbf{u}_{n}^{\prime\prime}\rightarrow\mathbf{u}_{0}.$ Then
\[
d(\mathbf{u}_{n}^{\prime\prime},\Gamma)=d(g_{n}\mathbf{u}_{n}^{\prime\prime
},\Gamma)\leq d(g_{n}\mathbf{u}_{n}^{\prime\prime},\mathbf{u}_{0}%
)\rightarrow0
\]
and this contradicts (\ref{kaka}).

Now we prove the other implication $\Leftarrow.$ Let $\mathbf{u}_{n}$ be a
sequence such that $d(\mathbf{u}_{n},\Gamma)\rightarrow0,$ then there exists
$\mathbf{v}_{n}\in\Gamma$ s.t.
\begin{equation}
d(\mathbf{u}_{n},\Gamma)\geq d(\mathbf{u}_{n},\mathbf{v}_{n})-\frac{1}{n}.
\label{triplo}%
\end{equation}

Since $V$ is G-compact, also $\Gamma$ is G-compact; so, for a suitable
sequence $g_{n}$, we have $g_{n}\mathbf{v}_{n}\rightarrow\mathbf{\bar{w}}%
\in\Gamma.$ We get the conclusion if we show that $V(\mathbf{u}_{n}%
)\rightarrow0.$ We have by (\ref{triplo}), that $d(\mathbf{u}_{n}%
,\mathbf{v}_{n})\rightarrow0$ and hence $d(g_{n}\mathbf{u}_{n},g_{n}%
\mathbf{v}_{n})\rightarrow0$ and so, since $g_{n}\mathbf{v}_{n}\rightarrow
\mathbf{\bar{w},}$ we have $g_{n}\mathbf{u}_{n}\rightarrow\mathbf{\bar{w}}%
\in\Gamma.$ Therefore, by the continuity of $V$ and since $\mathbf{\bar{w}}%
\in\Gamma,$ we have $V\left(  g_{n}\mathbf{u}_{n}\right)  \rightarrow V\left(
\mathbf{\bar{w}}\right)  =0$ and we can conclude that $V\left(  \mathbf{u}%
_{n}\right)  \rightarrow0.$

$\square$

In the cases in which we are interested, $X$ is an infinite dimensional
manifold; then if you choose generic $e_{0}$ and $c_{0},$ $V$ is not
$G$-compact since the set $\Gamma=\left\{  \mathbf{u}\in X:E(\mathbf{u}%
)=e_{0},\ C(\mathbf{u})=c_{0}\right\}  $ has codimension 2. However, Th.
(\ref{furbo}) allows to determine $e_{0}$ and $c_{0}$ in such a way that $V$
is $G $-compact and hence to prove the existence of solitons by using Theorem
\ref{astraco}.

\textbf{Proof of Th. \ref{astra1}.} In order to prove Th. \ref{astra1}, we
will use Th. \ref{astraco} with $e_{0}=e_{\delta}$ and $c_{0}=c_{\delta}$
where $e_{\delta}$ and $c_{\delta}$ are given by Th. \ref{furbo}.

We set
\begin{equation}
V\left(  \mathbf{u}\right)  =\left(  E\left(  \mathbf{u}\right)  -e_{\delta
}\right)  ^{2}+\left(  C\left(  \mathbf{u}\right)  -c_{\delta}\right)  ^{2}.
\label{liap}%
\end{equation}
We show that $V$ is $G$-compact: let $\mathbf{w}_{n}$ be a minimizing sequence
for $V,$ then $V\left(  \mathbf{w}_{n}\right)  \rightarrow0$ and consequently
$E\left(  \mathbf{w}_{n}\right)  \rightarrow e_{\delta}$ and $C\left(
\mathbf{w}_{n}\right)  \rightarrow c_{\delta}$. Now, since
\[
\min J_{\delta}=\frac{e_{\delta}}{c_{\delta}}+\delta\left[  e_{\delta
}+ac_{\delta}^{s}\right]  ,
\]
we have that $\mathbf{w}_{n}$ is a minimizing sequence also for $J_{\delta}. $
Then, since $J_{\delta}$ is $G$-compact, we get
\begin{equation}
\mathbf{w}_{n}\ \text{is}\ G\text{-compact}. \label{paracula}%
\end{equation}
So we conclude that $V$ is $G$-compact and hence the conclusion follows by
using Theorem \ref{astraco}.

$\square$

\bigskip

\section{Some applications\label{application}}

In this section we show how the abstract theorem \ref{astra1} can be used to
prove existence of solitons in some dynamical systems, namely for the
nonlinear Schr\"{o}dinger equation, the nonlinear wave equation, the nonlinear
beam equation and the Klein-Gordon-Maxwell equation. We do not carry out
detailed proofs since they are contained in ( \cite{BBBM},\ \cite{BBGM},
\cite{milano}, \cite{befobeam}, \cite{befolat}, \cite{befoQ}).

\subsection{The nonlinear Schroedinger equation}

Let us consider the nonlinear Schroedinger equation (NSE):
\begin{equation}
i\partial_{t}\psi=-\frac{1}{2}\Delta\psi+\frac{1}{2}W^{\prime}(\psi)
\tag{NSE}\label{NSE}%
\end{equation}

The solutions of this equations are critical points of the functional%
\[
\mathcal{S}=\int\left[  \operatorname{Re}(i\partial_{t}\psi\bar{\psi}%
)-\frac{1}{2}\left\vert \nabla\psi\right\vert ^{2}-W(\psi)\right]  dx\ dt
\]
Noether's theorem states that any invariance for a one-parameter group of the
Lagrangian implies the existence of an integral of motion (see e.g.
\cite{Gelfand} or \cite{milano}). They are derived by a continuity equation.

Now we describe the first integrals which will be relevant for this paper,
namely the energy, the \textquotedblleft hylenic charge\textquotedblright\ and
the momentum.

\textbf{Energy.}\emph{\ }The energy, by definition, is the quantity which is
preserved by the time invariance of the Lagrangian; it has the following form
\begin{equation}
E(\psi)=\int\left[  \frac{1}{2m}\left\vert \nabla\psi\right\vert ^{2}%
+W(\psi)\right]  dx.
\end{equation}

\textbf{Hylenic charge. }Following \cite{milano} the \textit{hylenic charge},
is defined as the quantity which is preserved by the invariance of the
Lagrangian with respect to the action
\[
\psi\mapsto e^{i\theta}\psi.
\]
Thus, in this case, the charge is nothing else but the $L^{2}$ norm, namely:
\[
C(\psi)=\int\left\vert \psi\right\vert ^{2}dx
\]

The phase space is given by%
\[
X=H^{1}(\mathbb{R}^{N},\mathbb{C})
\]
and the generic point in $X$ will be denoted by%
\[
\mathbf{u}=\psi
\]
The norm of $X$ is given by%
\[
\left\Vert \mathbf{u}\right\Vert =\left(  \int\left(  \left\vert \nabla
\psi\right\vert +\left\vert \psi\right\vert ^{2}\right)  dx\right)  ^{\frac
{1}{2}}%
\]
and the auxiliary seminorm is given by
\[
\left\Vert \mathbf{u}\right\Vert _{\natural}=\ \underset{z\in\mathbb{R}^{N}%
}{\sup}\left(  \int_{B_{1}(z)}\left\vert \psi\right\vert ^{2}dx\right)
^{\frac{1}{2}}%
\]
Applying theorem \ref{astra1} we get the following existence result (see e.g.
\cite{CL82},\ \cite{BBGM}):

\begin{theorem}
\textsc{\ }Assume that%
\[
W(s)=\frac{1}{2}as^{2}+N(s)
\]
and that%
\begin{equation}
|N^{\prime}(s)|\leq c_{1}|s|^{q-1}+c_{2}|s|^{p-1}\text{ for some }2<q\leq
p<2^{\ast}. \label{Fp}%
\end{equation}%
\begin{equation}
N(s)\geq-c_{1}s^{2}-c_{2}|s|^{\gamma}\text{ for some }c_{1},c_{2}%
\geq0,\ {\ \gamma<2+\frac{4}{N}} \label{F0}%
\end{equation}%
\[
N(s_{0})<0\text{ for some }s_{0}>0
\]
Then NSE has a continuous family of hylomorphic solitons.
\end{theorem}

We can see that the assumptions of Th. \ref{astra1} but (EC-3) are easy to
verify. In order to verify (EC-3) we need the Nash inequality in the following
form
\begin{equation}
||\psi||_{L^{p}}^{p}\leq b_{p}||\psi||_{L^{2}}^{r}||\nabla\psi||_{L^{2}}^{q}
\label{nash}%
\end{equation}

where%
\begin{align*}
p  &  <2+\frac{4}{N},\\
q  &  =pN\left(  \frac{1}{2}-\frac{1}{p}\right) \\
r  &  =p-q.
\end{align*}

then, for $s=p/2$ and $a$ sufficiently large, we have that
\begin{align*}
E(\mathbf{u})+a\left\vert C(\mathbf{u})\right\vert ^{s}  &  \geq\frac{1}%
{2}||\nabla\psi||_{L^{2}}^{2}+\int W\left(  |\psi|\right)  +a||\psi||_{L^{2}%
}^{2s}\geq(\text{by (\ref{F0})}\\
&  \geq\frac{1}{2}||\nabla\psi||_{L^{2}}^{2}-c_{3}||\psi||_{L^{p}}^{p}%
+a||\psi||_{L^{2}}^{2s}\geq(\text{by (\ref{nash})}\geq0
\end{align*}

\bigskip

\subsection{The nonlinear wave equation}

The nonlinear wave equation (NWE) has the following structure
\begin{equation}
\square\psi+W^{\prime}\left(  \psi\right)  =0\tag{NWE}\label{NWE}%
\end{equation}
where $\psi:\mathbb{R}^{4}\mathbb{\rightarrow C}$ and where, with some abuse
of ntation we have set
\[
W^{\prime}(\psi)=F^{\prime}(\left\vert \psi\right\vert )\frac{\psi}{\left\vert
\psi\right\vert }%
\]
for some smooth function $F:\left[  0,\infty\right)  \rightarrow\mathbb{R}.$
The NWE has a variational structure, namely it is the Euler-Lagrange equation
with respect to the functional%
\[
S=\frac{1}{2}\int\int\left(  |\partial_{t}\psi|^{2}-|\nabla\psi|^{2}\right)
dx\ dt-\int\int W(\psi)dx\ dt
\]

The energy and the hylenic charge take the following form
\[
E=\frac{1}{2}\int\left(  |\partial_{t}\psi|^{2}+|\nabla\psi|^{2}\right)
dx+\int W(\psi)dx
\]%
\[
C=\operatorname{Im}\int\psi_{t}\overline{\psi}\ dx.
\]

In order to describe the phase space $X$, we rewrite NWE as an Hamiltonian
system:%
\[
\left\{
\begin{array}
[c]{c}%
\partial_{t}\psi=\phi\\
\\
\partial_{t}\phi=\Delta\psi-W^{\prime}(\psi);
\end{array}
\right.
\]
so we can see that the phase space is given by%
\[
X=H^{1}(\mathbb{R}^{N},\mathbb{C})\times L^{2}(\mathbb{R}^{N},\mathbb{C})
\]
and the generic point in $X$ will be denoted by%
\[
\mathbf{u}=\left[
\begin{array}
[c]{c}%
\psi\\
\phi
\end{array}
\right]
\]
The norm of $X$ is given by%
\[
\left\Vert \mathbf{u}\right\Vert =\left(  \int\left(  |\phi|^{2}+|\nabla
\psi|^{2}+|\psi|^{2}\right)  dx\right)  ^{\frac{1}{2}}%
\]
the auxiliary seminorm is given by%
\[
\left\Vert \mathbf{u}\right\Vert _{\natural}=\ \underset{z\in\mathbb{R}^{N}%
}{\sup}\left(  \int_{B_{1}(z)}|\psi|^{2}\ dx\right)  ^{\frac{1}{2}}%
\]

The following existence theorem holds (see \cite{BBBM})

\begin{theorem}
Assume that%
\[
W(s)=\frac{1}{2}m^{2}s^{2}+N(s)
\]
and that

\begin{itemize}
\item (W-i) \textbf{(Positivity}) $W(s)\geq0$

\item (W-ii) \textbf{(Nondegeneracy}) $W^{\prime\prime}(0)=m^{2}$\ $>0$

\item (W-iii) \textbf{(Hylomorphy}) $\exists s_{0}:\;N(s_{0})<0$\ 

\item (W--iiii) \textbf{(Growth}) there is a constant $c>0$ such that
\[
N^{\prime}(s)\geq-c_{1}s-c_{2}s^{p-1},\ \ 2<p<2^{\ast}%
\]
Then NWE has a continuous family of hylomorphic solitons.
\end{itemize}
\end{theorem}

\bigskip

\subsection{The nonlinear beam equation}

Let us consider the nonlinear beam equation (NBE)%

\begin{equation}
\frac{\partial^{2}u}{\partial t^{2}}+\frac{\partial^{4}u}{\partial x^{4}%
}+W^{\prime}(u)=0. \tag{NBE}\label{NBE}%
\end{equation}
where $u:\mathbb{R}^{2}\rightarrow\mathbb{R},$ $W:\mathbb{R}\rightarrow
\mathbb{R}.$

Equation (\ref{NBE}), with a suitable choise of $W(s),$ has been proposed as
model for a suspension bridge (see \cite{mkw87}, \cite{lam}, \cite{lambis},
\cite{McKW90}).

NBE has a variational structure, namely it is the Euler-Lagrange equation with
respect to the functional%
\[
S=\frac{1}{2}\int\int\left(  u_{t}^{2}-u_{xx}^{2}\right)  dx\ dt-\int\int
W(u)dx\ dt.
\]

The energy and the momentum take the following form
\[
E=\frac{1}{2}\int\left(  u_{t}^{2}+u_{xx}^{2}\right)  dx+\int W(u)dx
\]%
\[
C=-\int u_{t}u_{x}\ dx.
\]
In this case the momentum will play the role of hylenic charge. Equation
(\ref{NBE}), can be rewritten as an Hamiltonian system as follows:%
\[
\left\{
\begin{array}
[c]{c}%
\partial_{t}u=v\\
\\
\partial_{t}v=-\partial_{x}^{4}u-W^{\prime}(u)
\end{array}
\right.
\]

The phase space is given by%
\[
X=H^{2}(\mathbb{R}^{N},\mathbb{C})\times L^{2}(\mathbb{R}^{N},\mathbb{C})
\]
and the generic point in $X$ will be denoted by%
\[
\mathbf{u}=\left[
\begin{array}
[c]{c}%
u\\
v
\end{array}
\right]  .
\]

The norm of $X$ is given by%
\[
\left\Vert \mathbf{u}\right\Vert =\left(  \int\left(  v^{2}+u_{xx}^{2}%
+u^{2}\right)  dx\right)  ^{\frac{1}{2}}.
\]
The auxiliary seminorm is given by%
\[
\left\Vert \mathbf{u}\right\Vert _{\natural}=\left\Vert u\right\Vert
_{L^{\infty}}.
\]

The energy and the momentum (which in this case plays the role of hylenic
charge), as functionals defined on $X,$ take the following form%
\[
E\left(  \mathbf{u}\right)  =\frac{1}{2}\int\left(  v^{2}+u_{xx}^{2}\right)
dx+\int W(u)dx
\]%
\[
C\left(  \mathbf{u}\right)  =-\int vu_{x}\ dx.
\]
By using Theorem \ref{astra1} the following existence result can be proved
(\cite{befobeam})

\begin{theorem}
Assume that

\begin{itemize}
\item (W-i) \textbf{(Positivity}) $W(s)>0$ for $s\neq0$ and $\exists
\delta>0\ $such that $\left\vert s\right\vert \geq1\Rightarrow W(s)\geq\delta$

\item (W-ii) \textbf{(Nondegeneracy}) $W^{\prime\prime}(0)$\ $>0$

\item (W-iii) \textbf{(Hylomorphy}) $\exists M>0,\ \exists\alpha\in
\lbrack0,2),\forall s\geq0,$
\[
W(s)\leq M\left\vert s\right\vert ^{\alpha}%
\]
Then NWE has a continuous family of hylomorphic solitons.
\end{itemize}
\end{theorem}

\subsection{The Klein-Gordon-Maxwell equations}

Now let us consider the Klein-Gordon-Maxwell equations:
\begin{equation}
D_{\varphi}^{2}\psi-D_{\mathbf{A}}^{2}\psi+W^{\prime}(\psi)=0 \label{e1}%
\end{equation}%
\begin{equation}
\nabla\cdot\left(  \partial_{t}\mathbf{A}+\nabla\varphi\right)
=q\operatorname{Re}\left(  iD_{\varphi}\psi\overline{\psi}\right)  \label{e2}%
\end{equation}%
\begin{equation}
\nabla\times\left(  \nabla\times\mathbf{A}\right)  +\partial_{t}\left(
\partial_{t}\mathbf{A}+\nabla\varphi\right)  =q\operatorname{Re}\left(
iD_{\mathbf{A}}\psi\overline{\psi}\right)  . \label{e3}%
\end{equation}
Here $q$ denotes a positive parameter which, in some physical models,
represents the unit electric charge, $\nabla\times$ and $\nabla$ denote
respectively the curl and the gradient operators;
\[
\mathbf{%
A%
}=\mathbf{%
\mathbf{(}%
}A_{1},A_{2},A_{3}\mathbf{%
)%
}\in\mathbb{R}^{3}\text{ and }\varphi\in\mathbb{R}%
\]
are the gauge potentials;
\[
D_{\varphi}\psi=\left(  \partial_{t}+iq\varphi\right)  \psi
\]
is the covariant derivative with respect to the $t$ variable, and
\[
D_{\mathbf{A}}\psi=\left(  \nabla-iq\mathbf{A}\right)  \psi
\]
is the covariant derivative with respect to the $x$ variable (see for example
\cite{befogranas} and \cite{yang}).

Equations (\ref{e1}),.., (\ref{e3}) are invariant with respect to the gauge
transformations%
\begin{equation}
\psi\rightarrow e^{iq\chi}\psi\label{ga}%
\end{equation}%
\begin{equation}
\varphi\rightarrow\varphi-\partial_{t}\chi\label{ge}%
\end{equation}%
\begin{equation}
\mathbf{A\rightarrow A}+\nabla\chi\label{gi}%
\end{equation}
where $\chi\in C^{2}\left(  \mathbb{R}^{4}\right)  $.

If we make the following change of variables:
\begin{equation}
\mathbf{E=-}\left(  \partial_{t}\mathbf{%
A%
}+\nabla\varphi\right)  \label{pos1}%
\end{equation}%
\begin{equation}
\mathbf{H}=\nabla\times\mathbf{A}\label{pos2}%
\end{equation}%
\begin{equation}
\rho=-q\operatorname{Re}\left(  iD_{\varphi}\psi\overline{\psi}\right)
\label{caricona}%
\end{equation}%
\begin{equation}
\mathbf{j}=q\operatorname{Re}\left(  iD_{\mathbf{A}}\psi\overline{\psi
}\right)  .
\end{equation}
you get the following equations:%

\[
\nabla\cdot\mathbf{E}=\rho,\ \ \mathbf{Gauss\ equation}%
\]%
\[
\nabla\times\mathbf{H}-\frac{\partial\mathbf{E}}{\partial t}=\mathbf{j}%
,\ \ \mathbf{Amp\acute{e}re\ equation}%
\]%
\[
\nabla\times\mathbf{E+}\frac{\partial\mathbf{H}}{\partial t}%
=0,\ \ \mathbf{Faraday\ equation}%
\]%
\[
\nabla\cdot\mathbf{H}=0,\ \ \mathbf{No}\text{-}\mathbf{monopole\ equation}%
\]

We write $\psi$ in polar form
\begin{equation}
\psi(x,t)=u(x,t)\,e^{iS(x,t)},\;\;u\geq0,\;\;S\in\mathbb{R}/2\pi\mathbb{Z.}
\label{giulia}%
\end{equation}

\bigskip

Equation (\ref{e1}) can be split in the two following ones
\begin{equation}
\square u+W^{\prime}(u)+\left[  \left\vert \nabla S-q\mathbf{A}\right\vert
^{2}-\left(  \partial_{t}S+q\phi\right)  ^{2}\right]  \,u=0 \label{e1+}%
\end{equation}%
\begin{equation}
\frac{\partial}{\partial t}\left[  \left(  \partial_{t}S+q\phi\right)
u^{2}\right]  -\nabla\cdot\left[  \left(  \nabla S-q\mathbf{A}\right)
u^{2}\right]  =0. \label{e1cont}%
\end{equation}
and, using the variables $\mathbf{j}$ and $\rho,$ these equations can be
written as follows:%
\[
\square u+W^{\prime}(u)+\frac{\mathbf{j}^{2}-\rho^{2}}{q^{2}u^{3}}=0
\]%
\[
\partial_{t}\rho+\nabla\cdot\mathbf{j}=0.
\]

Thus, our equations, in the gauge invariant variables, become:
\[
\square u+W^{\prime}(u)+\frac{\mathbf{j}^{2}-\rho^{2}}{q^{2}u}%
=0,\ \ \mathbf{Matter\ equation}%
\]%
\[
\nabla\cdot\mathbf{E}=\rho,\ \ \mathbf{Gauss\ equation}%
\]%
\[
\nabla\times\mathbf{H}-\frac{\partial\mathbf{E}}{\partial t}=\mathbf{j}%
,\ \ \mathbf{Amp\acute{e}re\ equation}%
\]%
\[
\nabla\times\mathbf{E+}\frac{\partial\mathbf{H}}{\partial t}%
=0,\ \ \mathbf{Faraday\ equation}%
\]%
\[
\nabla\cdot\mathbf{H}=0,\ \ \mathbf{No}\text{-}\mathbf{monopole\ equation}%
\]

Peculiar difficulties with the NKGM equations:

\begin{itemize}
\item (j) the Lagrangian density
\begin{align*}
\mathcal{L(}\psi,\partial_{t}\psi,\mathbf{A,}\partial_{t}\mathbf{A,}\varphi)
&  =\frac{1}{2}\left(  \left\vert D_{\varphi}\psi\right\vert ^{2}-\left\vert
D_{\mathbf{A}}\psi\right\vert ^{2}\right)  -W(\psi)\\
&  +\frac{1}{2}\left(  \left\vert \partial_{t}\mathbf{A+}\nabla\varphi
\right\vert ^{2}-\left\vert \nabla\times\mathbf{A}\right\vert ^{2}\right)  .
\end{align*}
is incomplete, namely it does not depend explicitely on $\partial_{t}\varphi$

\item (jj) the energy
\[
E\left(  \mathbf{u}\right)  =\frac{1}{2}\int\left(  \left\vert \partial
_{t}u\right\vert ^{2}+\left\vert \nabla u\right\vert ^{2}+\frac{\rho
^{2}+\mathbf{j}^{2}}{q^{2}u^{2}}+\mathbf{E}^{2}+\mathbf{H}^{2}\right)  dx+\int
W(u)dx
\]
does not have the "usual" form
\[
\left\{  \text{energy}\right\}  =\left\{  \text{positive quadratic
form}\right\}  +\left\{  \text{higher order terms}\right\}  .
\]

\end{itemize}

To overcome the difficulty (j), we set%

\begin{equation}
M_{L}=\left\{  (\mathbf{Q},\partial_{t}\mathbf{Q}):\partial_{t}\varphi
+\nabla\cdot\mathbf{A}=0\mathbf{,\ }\nabla\cdot\left(  \partial_{t}%
\mathbf{A+}\nabla\varphi\right)  =q\operatorname{Re}\left(  iD_{\varphi}%
\psi\overline{\psi}\right)  \right\}  \label{m0}%
\end{equation}
where%
\[
\mathbf{Q}=\left(  \psi\mathbf{,A},\varphi\right)  ,\ \partial_{t}%
\mathbf{Q}=\left(  \partial_{t}\psi,\partial_{t}\mathbf{A,}\partial_{t}%
\varphi\right)
\]
and we consider the \textit{modified Lagrangian}%
\begin{align}
\mathcal{\hat{L}}\mathbf{(Q,}\partial_{t}\mathbf{\mathbf{Q})}  &  =\frac{1}%
{2}\left[  \left\vert D_{\varphi}\psi\right\vert ^{2}-\left\vert
D_{\mathbf{A}}\psi\right\vert ^{2}\right]  -W(\psi)\nonumber\\
&  +\frac{1}{2}\left[  \left\vert \partial_{t}\mathbf{A}\right\vert
^{2}-\left\vert \nabla\mathbf{A}\right\vert ^{2}-\left(  \partial_{t}%
\varphi\right)  ^{2}+\left\vert \nabla\varphi\right\vert ^{2}\right]  .
\end{align}

The dynamics induced by $\mathcal{\hat{L}}$ is given by the following
equations:%
\begin{align}
D_{\varphi}^{2}\psi-D_{\mathbf{A}}^{2}\psi+W^{\prime}(\psi)  &
=0\label{dalembert1}\\
\square\mathbf{A}-q\operatorname{Re}\left(  iD_{\mathbf{A}}\psi\overline{\psi
}\right)   &  =0\label{dalembert2}\\
\square\varphi+q\operatorname{Re}\left(  iD_{\varphi}\psi\overline{\psi
}\right)   &  =0 \label{dalembert3}%
\end{align}
The following result is well known (for a proof see e.g. \cite{prodi})

\begin{theorem}
The set $M_{L}$ is invariant for the dynamics induced by equations
(\ref{dalembert1},\ref{dalembert2},\ref{dalembert3}). Moreover, if the initial
data are in $M_{L}$ and $\left(  \psi\mathbf{,A},\varphi\right)  $ is a smooth
solution of eq. (\ref{dalembert1},\ref{dalembert2},\ref{dalembert3}) then it
is also a solution of NKGM.
\end{theorem}

Now it is possible to define the conjugate variable of $\mathbf{Q}$ via the
Legendre transform. These conjugate variables will be denoted by
$\mathbf{P}=\left(  \hat{\psi},\widehat{\mathbf{A}},\hat{\varphi}\right)  .$
We have:%

\begin{equation}
\hat{\psi}=\frac{\partial\widehat{\mathcal{L}}}{\partial\left(  \partial
_{t}\psi\right)  }=\partial_{t}\psi+iq\varphi\psi=D_{\varphi}\psi\label{acan}%
\end{equation}

\begin{equation}
\widehat{\mathbf{A}}=\frac{\partial\widehat{\mathcal{L}}}{\partial\left(
\partial_{t}\mathbf{A}\right)  }=\partial_{t}\mathbf{A} \label{Ecan}%
\end{equation}%
\begin{equation}
\hat{\varphi}=\frac{\partial\widehat{\mathcal{L}}}{\partial\left(
\partial_{t}\varphi\right)  }=-\partial_{t}\varphi. \label{ican}%
\end{equation}
Now, the Hamiltonian is well defined and takes the form:%
\begin{equation}
\mathcal{H}\mathbf{(Q,P)}=\frac{1}{2}\int\left(  \left\vert \hat{\psi
}\right\vert ^{2}+\left\vert D_{\mathbf{A}}\psi\right\vert ^{2}+\left\vert
\nabla\varphi+\widehat{\mathbf{A}}\right\vert ^{2}+\left\vert \nabla
\times\mathbf{A}\right\vert ^{2}\right)  +\int W(\psi)
\end{equation}

Now let us see how to overcome the difficulty (jj). Using the gauge
independent variables, the energy takes the form
\[
E(\mathbf{u)}=\frac{1}{2}\int\left[  v^{2}+\left\vert \nabla u\right\vert
^{2}+\frac{\rho^{2}+\mathbf{j}^{2}}{q^{2}u^{2}}+\mathbf{E}^{2}+\mathbf{H}%
^{2}\right]  dx+\int W(u)dx.
\]
where $v=\partial_{t}u.$

We introduce new gauge invariant variables which eliminate this singularity:%
\begin{equation}
\theta=\frac{-\rho}{qu};\ \Theta=\frac{\mathbf{j}}{qu}. \label{theta}%
\end{equation}
Using these new variables the energy and the charge take the form:%
\begin{equation}
E\left(  \mathbf{u}\right)  =\frac{1}{2}\int\left[  v^{2}+\left\vert \nabla
u\right\vert ^{2}+\theta^{2}+\Theta^{2}+\mathbf{E}^{2}+\mathbf{H}^{2}\right]
+\int W(u). \label{E1}%
\end{equation}

\begin{equation}
C\left(  \mathbf{u}\right)  =-q\int\theta udx. \label{C1}%
\end{equation}

Now we can define the functional framework which allows to apply Th.
\ref{astra1}

We can take a norm defined by the quadratic part of the energy, namely%
\begin{equation}
\left\Vert \mathbf{u}\right\Vert =\left(  \int\left[  v^{2}+\left\vert \nabla
u\right\vert ^{2}+m^{2}u^{2}+\theta^{2}+\Theta^{2}+\mathbf{E}^{2}%
+\mathbf{H}^{2}\right]  \ dx\right)  ^{\frac{1}{2}}. \label{norm}%
\end{equation}

and the auxiliary seminorm:%
\[
\left\Vert \mathbf{u}\right\Vert _{\sharp}=\ \underset{z\in\mathbb{R}^{3}%
}{\sup}\left(  \int_{B_{1}(z)}u^{2}\ dx\right)  ^{\frac{1}{2}}.
\]

So we get  the following functional framework. We will denote by $V$ the
completion of $C_{0}^{\infty}(\mathbb{R}^{3},\mathbb{R}^{12})$ with respect to
the norm (\ref{norm}) so that
\[
\mathbf{u}=\left(  u,v,\theta,\Theta,\mathbf{E},\mathbf{H}\right)  \in V\cong
H^{1}(\mathbb{R}^{3})\times L^{2}(\mathbb{R}^{3},\mathbb{R}^{11})
\]
Finally $X\subset V$ will denote the closure of
\[
X_{0}=\left\{  \mathbf{u}\in C^{\infty}(\mathbb{R}^{3},\mathbb{R}%
^{12})\mathbf{\ |\ }\nabla\cdot\mathbf{E}=\rho,\ \nabla\cdot\mathbf{H}%
=0,\ E(\mathbf{u})<+\infty\right\}
\]
with respect to the norm $\left\Vert \mathbf{u}\right\Vert $.

By using Theorem \ref{astra1}, the following existence result can be proved
(\cite{befoQ})

\begin{theorem}
Assume that%
\[
W(s)=\frac{1}{2}m^{2}s^{2}+N(s)
\]
and that

\begin{itemize}
\item (W-i) \textbf{(Positivity}) $W(s)\geq0$

\item (W-ii) \textbf{(Nondegeneracy}) $m^{2}$\ $>0$

\item (W-iii) \textbf{(Hylomorphy}) $\exists\bar{s}>0\ $and $\alpha\in\left(
0,m\right)  $ such that $W(\bar{s})\leq\frac{1}{2}\alpha^{2}\bar{s}^{2}$

\item (W-iiii)\textbf{(Growth condition}) There are constants $a,b>0,$ $6>p>2
$ s.t. $|N^{\prime}(s)|\ \leq as^{p-1}+bs^{2-\frac{2}{p}}$.
\end{itemize}

Then NKGM have a continuous family of hylomorphic solitons.
\end{theorem}

\end{document}